\documentclass{amsart}
\usepackage{amsmath,amsthm,amsfonts,amssymb,amscd}
\usepackage[latin1]{inputenc}
\usepackage{graphics}
\usepackage{graphicx}

\begin{document}

\newcommand{\la}{{\lambda}}
\newcommand{\po}{{\partial}}
\newcommand{\ov}{\overline}
\newcommand{\re}{{\mathbb{R}}}
\newcommand{\nb}{{\mathbb{N}}}
\newcommand{\lc}{{\mathcal{L}}}
\newcommand{\Uc}{{\mathcal U}}
\newcommand{\fc}{{\mathcal F}}
\newcommand{\al}{{\alpha}}
\newcommand{\vr}{{\varphi}}
\newcommand{\om}{{\omega}}
\newcommand{\La}{{\Lambda}}
\newcommand{\be}{{\beta}}
\newcommand{\Om}{{\Omega}}
\newcommand{\ve}{{\varepsilon}}
\newcommand{\ga}{{\gamma}}
\newcommand{\Ga}{{\Gamma}}
\newcommand{\supp}{\operatorname{supp}}
\newcommand{\cb}{{\mathbb{C}}}
\newcommand{\gc}{{\mathcal G}}
\newcommand{\zb}{{\mathbb{Z}}}
\newcommand{\pc}{{\mathcal P}}
\newcommand{\oc}{{\mathcal O}}

\title{Positive Neighborhoods of Rational Curves}
\date{\today}

\author{M. Falla Luza}
\author{P. Sad}

\address{Instituto de Matem\' atica e Estat\'\i stica, Universidade Federal Fluminense, Rua M\' ario Santos Braga s/n, 24020-140, Niter\'oi, Brazil}
\address{Instituto Nacional de Matematica Pura e Aplicada, IMPA, Estrada dona Castorina 110, Jardim Botânico, Rio de Janeiro, Brazil}
\address{IRMAR, Universit\'e de Rennes 1, 35042 Rennes Cedex, France}

\email{maycolfl@impa.br, sad@impa.br}

\begin{abstract} 
\noindent {We study neighborhoods of rational curves in surfaces with self-intersection number 1 that can be 
linearised}.
\end{abstract}

\maketitle

\section{Introduction}

Let us consider a holomorphic embedding of the complex projective line into a complex surface. The first invariant we can attach to this 
embedding is the self-intersection number of the image. When $C^2>0$, there is a big moduli space of germs of neighborhoods $(S,C)$,
see \cite{Hur,Mish}. In the paper \cite{FallaLoray} the authors give an extensive description of these
embeddings when the self-intersection number is equal to 1. In particular they describe completely the embeddings that are equivalent to the
one of the complex projective line as a line of the projective plane. Our interest here is to present a geometric approach to this
case using holomorphic foliations, which perhaps may be useful to study other situations.

We take then a smooth rational curve $C$ embedded in a complex surface $S$ with self-intersection number $C\cdot C = 1$. Let $H$ be a line in the complex projective plane $\mathbb CP(2)$. One of the results of \cite{FallaLoray} is the following:

\vspace {0.2in}
\noindent {\bf Theorem.}\, Suppose that there exist three different holomorphic fibrations over $C$ in $S$. Then a neighborhood of $C$ is holomorphically equivalent to a neighborhood of $H$ in $\mathbb CP(2)$.

\vspace {0.2in}
This is the result we study here using a different viewpoint.

We remark that any holomorphic fibration over $H$ in $\mathbb CP(2)$ is in fact a linear one, that is, defined by a pencil of lines of $\mathbb CP(2)$ with the base point lying outside $H$. In fact, a holomorphic fibration extends to a foliation on $\mathbb{P}^2$ which, being totally transverse to $H$, must be of degree $0$.

Our idea to prove the Theorem is to analyse the curve of tangencies between two fibrations over $C$. When this curve is a common fiber, we construct a holomorphic diffeomorphism from a neighborhood of $C$ to a neighborhood of $H$. The other possibilities are: (i)the curve of tangencies is $C$ itself; (ii)the curve of tangencies is transverse to $C$ but
it is not a common fiber. What we show is that in the presence of three different fibrations, necessarily two of them have a common fiber.

It should be mentioned that in \cite {FallaLoray} the authors describe also the possibility of existence of fibrations over rational curves of positive self-intersection number. We explain this point in Section 3.

We are grateful to Frank Loray for inumerous discussions around this subject.

\section{Curves of Tangencies}

For the reader's convenience, we present a short description of the curve of tangencies between two fibrations of $S$ over $C$.

Let us take a regular foliation $\mathcal F$ of the surface $S$ generically transverse to the (smooth) curve $C$. We will use the notation:

\begin{enumerate}
\item $N_{\mathcal F}$ is the normal bundle of $\mathcal F$.
\item $tang(\mathcal F,C)$ is the tangency locus between $\mathcal F$ and $C$.
\item If $\mathcal G$ is another foliation of $S$ generically transverse to $C$, then $tang(\mathcal F,\mathcal G)$ is the divisor of tangencies between $\mathcal F$ and $\mathcal G$.
\item $K_S$ is the canonical bundle of $S$.
\item $\chi(C)$ is the Euler characteristic of $C$.
\end {enumerate}

According to \cite{Brun}, we have:
$$
N_{\mathcal F}\cdot C= \chi(C)+ tang(\mathcal F,C) \,\,\,\,\,(N_{\mathcal G}\cdot C= \chi(C)+ tang(\mathcal G,C))
$$
and
$$
tang(\mathcal F, \mathcal G)\cdot C= N_{\mathcal F}\cdot C + N_{\mathcal G}\cdot C + K_{S}\cdot C
$$

We are interested in the case where $\mathcal F$ and $\mathcal G$ are fibrations over the rational curve $C$ which has self-intersection number equal to 1. From the formulae above we get then $N_{\mathcal F}\cdot C= N_{\mathcal G}\cdot C=2$. The adjunction formula gives
$$
-K_{S}\cdot C= \chi(C) + C\cdot C= 3
$$
We get then $tang(\mathcal F,\mathcal G)\cdot C = 1$. Therefore the possibilities are

\begin{itemize}
\item a (smooth, connected) curve of tangencies transversal to $C$
\item $tang(\mathcal F, \mathcal G) = C$
\end{itemize}

In both cases, the order of tangency between $\mathcal F$ and $\mathcal G$ is 1 along the curve of tangencies.
 
\section{The Case of a Common Fiber}

We consider now the first of the possibilities above with the curve of tangencies as a common fiber. We have then two
fibrations $\mathcal F$ and $\mathcal G$ over the curve $C$ with a common fiber $L$; we put $C \cap L= \{a\}$. We choose a line $H \subset {\mathbb C}P(2)$ and fix a biholomorphism $\phi: C \rightarrow H$. We take also two linear fibrations ${\mathcal F}^{\prime}$ and ${\mathcal G}^{\prime}$ in $\mathbb CP(2)$ which are transverse to $H$ and have a common fiber passing through $\phi(a)$.

\vspace {0.2in}
\noindent {\bf Proposition 1.}\, There exists a neighborhood of $C$ in $S$ which is biholomorphic to a neighborhood of $H$ in $\mathbb CP(2)$.

\begin {proof}

1) To each point $p\in S\setminus C$ \,close to $C$ we associate the points $\pi_{\mathcal F}(p)\in C$ and $\pi_{\mathcal G}(p)\in C$ as the points where the fibers of $\mathcal F$ and $\mathcal G$ passing through $p$ intersect $C$. We define $\Phi (p)\in \mathbb CP(2)$ as the point where the lines of ${\mathcal F}^{\prime}$ and ${\mathcal G}^{\prime}$ through $\phi(\pi_{\mathcal F}(p))$ and $\phi(\pi_{\mathcal G}(p))$ meet in $\mathbb CP(2)$ (later on we will need to make a choice for these pencils in $\mathbb CP(2)$). It is easy to see that $\Phi$ is a biholomorphism between neighborhoods of $S\setminus \{a\}$ and $H\setminus \phi(a)$. Let us show now that $\Phi$ extends along $L$ as a holomorphic map.

Let us take coordinates $(x,y)$ around the point $a\in C$ such that $(0,0)\in {\mathbb C}^2$ corresponds to $a$, $y=0$ parametrizes $C$ and $\mathcal F$ is defined as $d\,x=0$ (the common fiber $L$ becomes $x=0$). The fibration $\mathcal G$ is therefore defined by the set of curves
$$
x=x_0+A(x_0,y)
$$

\vspace {0,1in}
\noindent where $x_0\in \mathbb C$ is close to $0\in \mathbb C$, $A(x_0,y)$ is holomorphic and $A(0,y)\equiv 0$. In $\mathbb CP(2)$ we use affine coordinates $(X,Y)$ as to have the pencil ${\mathcal F}^{\prime}$ given by $d\,X=0$ and the pencil $\mathcal G$ given by the set of lines issued from the point $(0,1)$. Without any loss of generality we may
assume that $\phi(x,0)=(x,0)$. A simple computation shows that

$$
\Phi(x,y)=(x,1-\dfrac{x}{x_0})=(x,-\dfrac{A(x_0,y)}{x_0})
$$
where $(x_0,0)$ is the point of intersection of the $\mathcal G$-fiber through $(x,y)$ with the $x$-axis.

\vspace{0,1 in}
When $x \rightarrow 0$ (so $x_0 \rightarrow 0$), \,\,$\Phi(x,y)$ converges to $(0,-\dfrac{\delta A}{\delta x_0}(0,y))$.
Consequently, by Riemann´s extension theorem, $\Phi$ extends holomorphically to the fiber $L$. 

\vspace {0,1 in}
\noindent 2)In order to finish the proof of Proposition 1, we must prove that $\Phi$ is invertible. This is the place where we need to make a choice of the linear fibrations ${\mathcal F}^{\prime}$ and  ${\mathcal G}^{\prime}$. First of all we notice that the restrictions of tangent spaces $TS|_C$ and $T{\mathbb CP(2)}|_H$ are isomorphic, that is, for each $p\in C$ there exists a linear isomorphism $U_p$ from $TS|_C(p)$ to $T{\mathbb CP(2)}|_H(\phi(p))$ which is holomorhic in $p$. We define then in $\mathbb CP(2)$ the linear fibrations (over $H$) ${\mathcal F}^{\prime}$ and  ${\mathcal G}^{\prime}$ such that their fibers over $\phi(a)$ are $U_p(T{\mathcal F}_p)$ and $U_p(T{\mathcal G}_p)$ respectively. Finally we take a holomorphic diffeomorphism $\Psi$ from a neighborhood of $a\in C$ to a neighborhood  of $(0,0)\in \mathbb C^2$ such that $d\Psi|_{C}(p)=U_p$ for $p$ close to $a$. In other words, we have now local coordinates $(x,y)$ around $a\in C$ with the properties: i)$C$ is given by $y=0$; ii)$a$ corresponds to $(0,0)$ and iii) the fibrations $\mathcal F$ and $\mathcal G$ have tangent spaces along the $x$-axis (near $(0,0)$) equal to the linear fibers of ${\mathcal F}^{\prime}$ and  ${\mathcal G}^{\prime}$. 

In order to simplify the exposition, we will assume from now on that $\mathcal F$ is given by $d\,x=0$, ${\mathcal F}^{\prime}$ is given by $d\, X=0$, $\phi(x,0)=(x,0)$ and that ${\mathcal G}^{\prime}$ has $(0,1)$ as base point. We proceed then to compute the derivative of $\Phi(x,y)=(x,Y)$ at $(0,0)\in {\mathbb C}^2$ along the $y$-axis. Let $u(x)$ denote the slope of $\mathcal G$ (or ${\mathcal G}^{\prime}$) at the point $(x,0)$; $u$ is a holomorphic function with $u(0)=\infty$. Let $v(x,y)$ be the point in the $x$-axis which belongs to the $\mathcal G$-fiber through $(x,y)$. We have then
$$
\dfrac{y}{x-v(x,y)}\rightarrow u(v(x,y))
$$
when $y\rightarrow 0$ along the same fiber. 
  

\begin{figure}[t]
\begin{center}
\includegraphics[width=0.7\textwidth]{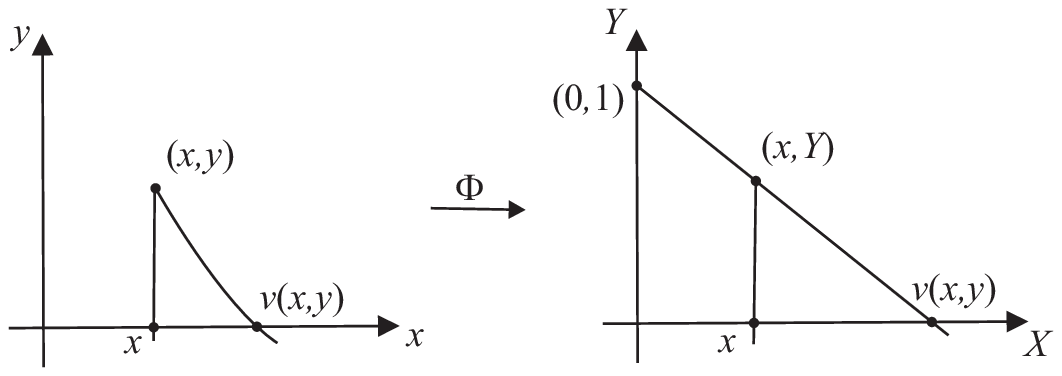}
\caption{}
\label{}
\end{center}
\end{figure}

\noindent But since 
$$
\dfrac{Y}{x-v(x,y)}=u(x,y)
$$ 
wee see that $\dfrac{Y}{y}\rightarrow 1$ when $y\rightarrow 0$ along the $\mathcal G$-fiber over $(x,y)$. When we make $x\rightarrow 0$,
we still have $\dfrac{Y}{y}\rightarrow 1$, which shows that $ d\Phi(0,0)$ restricted to the $y$-axis is the identity.
\end {proof}

Before closing this Section, let us turn to the situation (\cite {FallaLoray}) where the rational curve has positive self-intersection number.

\vspace {0.2in}
\noindent {\bf Proposition 2.}\, If $C^2 >1$ and there exist two fibrations $\mathcal F$ and $\mathcal G$  transverse to it then $C^2 = 2$ and there is a neighborhood of $C$ in $S$ which is equivalent to a neighborhood of the diagonal $\Delta$ in $\mathbb{P}^1 \times \mathbb{P}^1$. In particular, when  $C\cdot C=2$,  there are at most two fibrations transverse to it.

\begin{proof}
From the formula of Section $2$ we see that $tang(\mathcal F,\mathcal G)\cdot C= 2-C\cdot C$. We conclude that ${2- C\cdot C} \geq 0$, so that we have either $C\cdot C=1$ or $C\cdot C=2$. In this last case, we conclude also that the two fibrations are tranverse to each other and  to $C$ as well. We look now to the diagonal $\Delta$ inside ${\mathbb P}^1 \times {\mathbb P}^1$ (of course ${\Delta}\cdot {\Delta}=2$); let $\bar {\mathcal F}$ and $\bar {\mathcal G}$ be the horizontal and vertical fibrations to $\Delta$. We may construct as in the first part of the Proof of Proposition 1 the map $\Phi$ from  a neighborhood of $C$ to a neighborhood of $\Delta$ sending $C$ to $\Delta$, $\mathcal F$ to $\bar {\mathcal F}$ and $\mathcal G$ to $\bar {\mathcal G}$. This map is a diffeomorphism, since it is invertible. 
\end{proof}
\section{Disjoint Curves of Tangencies}

The first case we analyse in order to prove Theorem 1 is when we have two pairs of fibrations, say $(\mathcal F,\mathcal G)$ and $(\mathcal F,\mathcal H)$, with disjoint lines of tangencies $L_{\mathcal F,\mathcal G}$ and $L_{\mathcal F, \mathcal H}$ transverse to $C$ and to the respective pair of fibrations. Let us consider the map defined in the beginning of the proof of Proposition 1 (denoted here by $\Phi_{\mathcal F,\mathcal G}$). It conjugates the pair $(\mathcal F,\mathcal G)$ to a pair of linear pencils $({\mathcal F}^{\prime},{\mathcal G}^{\prime})$, which we take in affine coordinates as $d\,X=0$ and $(Y-1)d\,X-Xd\,Y=0$ in $\mathbb CP(2)$. We describe such a map in a neighborhood of $\{Q\}= L_{\mathcal F,\mathcal G} \cap C$.

We consider local coordinates $(x,y_1)$ such that: i)$(0,0)$ corresponds to $Q \in C$ and $y_1=0$ to $C$; ii) $\mathcal F$ is defined by $d\,x=0$. Let $y_1=l_{\mathcal F,\mathcal G}$ be the equation of $L_{\mathcal F, \mathcal G}$.

\vspace{0.2in}
\noindent {\bf Lemma 1.} There exists a local diffeomorphism $H(x,y_1)=(x,y)=(x,h(x,y_1))$ such that  $h(x,0)=0$ and transforms the pair $(d\,x=0,\,\mathcal G)$  to the pair
$$
(d\,x=0, \,\,\,d\,(x-y(y-2x))=0)
$$

\begin{proof}  We start assuming that $h(x,0)=0$; the fiber of $\mathcal G$ passing through each $(x,0)$ cuts the $y_1$-axis in two points $\{A^{+}(x),A^{-}(x)\}$ (with $A^{+}(0)=A^{-}(0)=0$), defining then an holomorphic involution $I_0$ that sends one point to the other (and $A(0)$ is fixed). The same is true for the fibration $d\,(x-y(y-2x))=0$ in relation to the $y$-axis: the fiber of $d\,(x-y(y-2x))=0$ through $(x,0)$ cuts the $y$-axis in the points $\{a^{+}(x), a^{-}(x)\}$ ($a^{+}(0)= a^{-}(0)=0$); the associated involution $i_0$ is conjugated to $I_0$ by  a holomorphic diffeomorphism $h_0(y_1)$ that satisfies $h_0(0)=0$.   

We may now extend holomorphically $h_0(y_1)$ to $h_{x}(y_1)$  in order to conjugate the involutions $I_x$ and $i_x$  that act in a close vertical (with $y_1=l_{\mathcal F, \mathcal G}(x)$ as the fixed point of each involution). 
\end{proof}

There is in fact noting special in the choice of $d\,(x-y(y-2x))=0$ to represent the fibration $\mathcal G$; it makes easier the analysis of $\Phi_{\mathcal F,\mathcal G}$. For example, the curve $L_{\mathcal F, \mathcal G}$ becomes the diagonal $y=x$.

As before we assume that the biholomorphism $\phi$ is $\phi(x,0)=(x,0)$. Then $\Phi_{\mathcal F,\mathcal G}$ is the rational map
$$
\Phi_{\mathcal F, \mathcal G}=(X,Y)=\left(x,1-\dfrac{x}{x-y(y-2x)}\right)=\left(x,\dfrac{-y(y-2x)}{x-y(y-2x)}\right)
$$
Let us consider in $\mathbb CP(2)$ the pencil $\mathcal P$ with center at some point $P\in H$, $P\neq Q$ and take its pull-back  ${\Phi_{\mathcal F, \mathcal G}}^*(\mathcal P)$; it is a foliation $\mathcal L$ defined in a neighborhood of $C$ with a radial singularity at the point $\Phi_{\mathcal F, \mathcal G}^{-1}(P)$ and a singularity at $Q$ of the form
$$
2x(y-x)d\,y - y^2d\,x + A(x,y)d\,y+ yB(x,y)d\,x =0
$$
where $A(x,y)=\sum_{i+j\geq 3}a_{ij}x^iy^j$ and $B(x,y)=\sum_{i+j\geq 2}b_{ij}x^iy^j$.
We remark the following fact that will be important for us: the curve of tangencies between $\mathcal L$ and a foliation defined in a neighborhood of $Q$ by $d\,y-y(1+C(x,y))d\,x=0$ ($C(x,y)=\sum_{i+j\geq 1}c_{ij}x^iy^j$) is, besides $y=0$, the curve 
$$
F(x,y)= (2x^2-A)(1+C)-B +y(1-2x(1+C))=0
$$   
We have then the intersection number at $Q$: $(F\cdot C)_Q\geq 2$.

\vspace{0.2in}
\noindent{\bf Proposition 3}.\,There are no fibrations $\mathcal F, \mathcal G$ and $\mathcal H$ such that  $L_{\mathcal F,\mathcal G}$ and $L_{\mathcal F, \mathcal H}$ are disjoint curves transverse to $C$ and to the respective fibrations. 
\begin {proof} We consider the pair of foliations $\mathcal L,\mathcal L^{\prime}$ constructed as we indicated above ($\mathcal L^{\prime}$ is obtained as the pull-back by $\Phi_{\mathcal F,\mathcal H}$ of the  pencil in $\mathbb CP(2)$  with center at the same point $P$; it has a singularity at the point $Q^{\prime}= L_{\mathcal F, \mathcal H}\cap C$). When we compute $tang(\mathcal L,\mathcal L^{\prime})\cdot C$, we have contributions at $Q$ and $Q^{\prime}$ (2 at least  for each point, as we remarked before) and at least 1 for the point $\Phi_{\mathcal F, \mathcal G}^{-1}(P)=\Phi_{\mathcal F, \mathcal H}^{-1}(P)$ (the $\mathcal F$-fiber through this point is a common separatrix of $\mathcal L,\mathcal L^{\prime}$); there is also the contribution of $C\cdot C=1$. All in all we get $tang(\mathcal L,\mathcal L^{\prime})\geq 6$.

But we have 
$$
tang(\mathcal L,\mathcal L^{\prime})\cdot C= N_{\mathcal L}\cdot C + N_{\mathcal L^{\prime}}\cdot C +K_S\cdot C
$$ 
From \cite[Chapter 2.2]{Brun} we have $N_{\mathcal L}\cdot C= Z(\mathcal L,C) + C\cdot C$, where $Z(\mathcal L,C)$ is the number of singularities of $\mathcal L$ in $C$ counted with multiplicities. It follows that $N_{\mathcal L}\cdot C=4$; in the same way $N_{{\mathcal L}^{\prime}}\cdot C=4$. Finally, since $K_S\cdot C=-3$, we arrive at $tang(\mathcal L,\mathcal L^{\prime})=5$, a contradiction.
\end{proof}

\noindent {\bf Remark 1}. Let us consider a pair of fibrations $(\mathcal F,\mathcal G)$ whose line of tangencies $L_{\mathcal F,\mathcal G}$ is transverse to $C$ at the point $Q$ and is not a common fiber. It may happen that  $L_{\mathcal F,\mathcal G}$ is tangent to the fibers of $\mathcal F$ and $\mathcal G$ at $Q$, but it is transverse to the fibers of $\mathcal F$ and $\mathcal G$ at any point close to but different from $Q$.

\section{Intersecting Curves of Tangencies}

We consider now the case where we have curves of tangencies $L_{\mathcal F ,\mathcal G}$ and $L_{\mathcal F, \mathcal H}$  transverse to $C$ and to the respective pair of fibrations, and have a common point $Q\in C$. We proceed as in Section 4 in order to define foliations $\mathcal L$ and ${\mathcal L}^{\prime}$ which leave $C$ invariant and have singularities at the point $Q$ and  radial singularities at another point of $C$. The same computation we did at the end of Proposition 1 informs us that $tang(\mathcal L,{\mathcal L}^{\prime})\cdot C= 5$.

We will show that a direct computation produces a different result. Let us then
analyse $tang(\mathcal L,{\mathcal L}^{\prime})\cdot C$ at the point $Q$. 

Let us consider a local coordinate $(x,y_1)$ around $Q$ such that $(0,0)$ corresponds to $Q$, $C$ becomes $y_1=0$ and $\mathcal F$ is written as $d\,x=0$; we still write $\mathcal G$ and $\mathcal H$ for the two other fibrations in these coordinates. From Lemma 1 we know that there exist two local diffeomorphisms $H(x,y_1),H^{\prime}(x,y_1)$ that transform  the pairs $(d\,x=0,\,\mathcal G)$ and $(d\,x=0,\,\mathcal H)$ to the pairs
$$
(d\,x=0, \,\,\,d\,(x-y(y-2x))=0)
$$
and
$$
(d\,x=0,\,\,\,d\,(x-y^{\prime}(y^{\prime}-2x)=0)
$$
respectively. 

\vspace{0.2in}
Now we take $\mathcal L$ and $\mathcal L^{\prime}$ as the pull-back´s by these maps $\Phi_{\mathcal F, \mathcal G}$ and $\Phi_{\mathcal F, \mathcal H}$ of the radial pencil $\mathcal P$ in $\mathbb P^2$ with center in some $P\in H$, $P\neq (0,0)$. The foliation $\mathcal L$ has a singularity at $Q$ given by  $2x(y-x)d\,y-y^2d\,x+h.o.t=0$ and another singularity at $\Phi_{\mathcal F,\mathcal G}^{-1}(P)$ of radial type; analogously for $\mathcal L^{\prime}$ in the coordinates $(x,y^{\prime})$.

Let us analyse the tangencies between $\mathcal L$ and $\mathcal L^{\prime}$ at the point $Q$, or the tangencies at $(0,0)$ between $2x(y-x)d\,y-y^2d\,x +h.o.t.=0$ and the pull-back of $2x(y^{\prime}-x)d\,y^{\prime}-{y^{\prime}}^2d\,x +h.o.t.=0$ by $L=H^{\prime}\circ H^{-1}$. We write $L(x,y)= (x,yC(x,y))$ and $L_t(x,y)=(x, C_t(x,y))=(x,y(C(x,y)+tx))$, for $t$ close to $0$; then we have $L_t(x,y)=(x,y( a+tx+\sum_{i\geq 1}a_i(x)y^i))$ with $a\neq0$. As before, we are interested in the components of the curve of tangencies different from $C$.

We now look for the tangencies between $\mathcal L: 2x(y-x)d\,y-y^2d\,x +h.o.t.=0$ and the pull-back ${\mathcal L^{\prime}}_t$ of $2x(y^{\prime}-x)d\,y^{\prime}-{y^{\prime}}^2d\,x +h.o.t.=0$ by $L_t$; we find a curve given by
$(4x^4 t^2+ \dots)y +y^2(\dots)=0$, so we get at least 4 when we compute $tang(\mathcal L,{\mathcal L^{\prime}}_t)\cdot C$ at $Q$ without taking $C$ into account. It follows that when we go to the limit $t=0$ we have also a contribution of at least 4 for $tang(\mathcal L,\mathcal L^{\prime})\cdot C$ at $Q$ without taking $C$ into account. Including the contributions of $C$ and of the singularities of $\mathcal L$ and $\mathcal L^{\prime}$ corresponding to $\Phi_{\mathcal F,\mathcal G}^{-1}(P)=\Phi_{\mathcal F,\mathcal H}^{-1}(P)$ (which have a common separatrix), we finally arrive at $tang(\mathcal L,\mathcal L^{\prime})\cdot C\geq 6$, a contradiction. Therefore we may state:

\vspace{0.2in}
\noindent {\bf Proposition 4.}\,\,There are no fibrations $\mathcal F, \mathcal G$ and $\mathcal H$ such that $L_{\mathcal F,\mathcal G}$ and $L_{\mathcal F, \mathcal H}$ are intersecting curves transverse to $C$ and to the respective pair of fibrations.

\section{Proof of the Theorem} 

A basic property of a embedded rational curve $C$ with self-intersection number equal to $1$ is that it can be deformed along a family of smooth rational curves of self-intersection numbers also equal to $1$ in some neighborhood of $C$, see \cite{Savelev} and also \cite{LB}.

We can now prove the Theorem by contradiction. Given three fibrations $\mathcal F$, $\mathcal G$ and $\mathcal H$ over $C$, we have two possibilities: i) either the three fibrations are tangent along $C$ or ii) there are two pairs which have lines of tangencies transverse to $C$  without being a common fiber of the respective pair (otherwise Proposition 1 applies). We then replace $C$ by a close smooth rational curve $C^{\prime}$ as explained above; the lines of tangencies will become transverse to the fibrations by Remark 1 (in the first case, the lines of tangencies coincide with $C$, of course). We are now in a setting where Propositions 3 and 4 may be applied to reach the contradiction.

\end{document}